\documentclass[final]{siamltex}

\usepackage{graphicx}
\usepackage{amssymb}
\usepackage{epstopdf}
\usepackage{cite}
\usepackage[hidelinks]{hyperref}
\usepackage{rotfig_tikz}
\usepackage{ 
  tikz, 
}
\usetikzlibrary{calc}
\usetikzlibrary{shapes}
\usetikzlibrary{fit}
\usetikzlibrary{arrows,decorations.markings,patterns}
\usepackage{ 
  pgfplots,
}

\pgfplotsset{%
  compat=1.8,%
}
\usepgfplotslibrary{statistics}

\definecolor{SPECorange}{rgb}{1.0,.5625,0}
\definecolor{SPECblue}{rgb}{0,0,0.75}
\definecolor{SPECred}{rgb}{0.75,0,0}
\definecolor{SPECgreen}{rgb}{0,0.75,0}
\definecolor{SPECblack}{rgb}{0.75,0.75,0.75}
\definecolor{SPECCorange}{rgb}{1,0.9,0.6}
\definecolor{SPECCblue}{rgb}{0,0,0.25}
\definecolor{SPECCred}{rgb}{1.0,0,0}
\definecolor{SPECCgreen}{rgb}{0.25,0.75,0.25}
\definecolor{SPECgray}{rgb}{0.35,0.35,0.35}

\usepackage{tabularx}
\usepackage{booktabs} 
\newcolumntype{C}[1]{>{\centering\arraybackslash}p{#1}} 
\usepackage{multirow}
\newcommand{\cplxs}{\mbox{$\mathbb{C}$}}

\newcommand{\cnxn}{\mbox{$\cplxs^{n \times n}$}}
\newcommand{\absval}[1]{\mbox{$\mid\!#1\!\mid$}}

\title{The RQR algorithm\thanks{This research was partially supported by the
    Research Council KU Leuven (Belgium), project C16/21/002 (Manifactor: Factor
    Analysis for Maps into Manifolds) and by the Fund for Scientific Research --
    Flanders (Belgium), projects G0A9923N (Low rank tensor approximation
    techniques for up- and downdating of massive online time series clustering)
    and G0B0123N (Short recurrence relations for rational Krylov and orthogonal
    rational functions inspired by modified moments).}}  \author{Daan
  Camps\footnotemark[3]\and Thomas Mach\footnotemark[4]\and Raf
  Vandebril\footnotemark[5]\and David~S.~Watkins\footnotemark[6]}

\begin{document}

\date{today}
\maketitle

\renewcommand{\thefootnote}{\fnsymbol{footnote}}

\footnotetext[3]{National Energy Research Scientific Computing Center, Lawrence Berkeley National Laboratory, USA (\texttt{dcamps@lbl.gov}).}%
\footnotetext[4]{University of Potsdam, Institute of Mathematics, 
  Karl-Liebknecht-Str.\ 24--25, 14476 Potsdam, Germany;
  \mbox{(\texttt{thomas.mach@uni-potsdam.de})}.}%
\footnotetext[5]{Department of Computer Science, KU Leuven,
  Belgium (\texttt{raf.vandebril@cs.kuleuven.be}).}
\footnotetext[6]{Department of Mathematics, Washington State University
   (\texttt{watkins@math.wsu.edu})}
\renewcommand{\thefootnote}{\arabic{footnote}}

\begin{abstract}
Pole-swapping algorithms, generalizations of bulge-chasing algorithms, have been
shown to be a viable alternative to the bulge-chasing QZ algorithm for solving the 
generalized eigenvalue problem for a matrix pencil $A - \lambda B$.   It is natural to try to
devise a pole-swapping algorithm that solves the standard eigenvalue problem for a 
single matrix $A$.   This paper introduces such an algorithm and shows that it
is competitive with Francis's bulge-chasing QR algorithm.
\end{abstract}



\section{Introduction}

The standard algorithm for computing the eigenvalues of a small to medium-sized
non-Hermitian matrix $A\in\cnxn$ is still Francis's implicitly-shifted QR
algorithm \cite{Fra61b,Wat11}.\footnote{LAPACK, the most widely used library for
  linear algebra computation, uses Francis's algorithm in
  \texttt{https://github.com/Reference-LAPACK/lapack/blob/master/SRC/zhseqr.f}.}
In many applications, eigenvalue problems arise naturally as generalized
eigenvalue problems for a pencil $A - \lambda B$, and for these problems the
Moler-Stewart variant of Francis's algorithm \cite{MolSte73}, commonly called
the QZ algorithm, can be used.  These are bulge-chasing algorithms~\cite{Wat07}:
they introduce a non-zero matrix element (or bulge) in an upper Hessenberg
matrix (resp.\ upper Hessenberg-triangular pencil), and chase it downwards along
the subdiagonal.

In recent years a generalization of QZ called RQZ was introduced by Camps,
Meerbergen, and Vandebril \cite{CaMeVa19a,Cam19}.  See also Steel et al.\
\cite{StCaMeVa21} and Camps et al.\ \cite{CaMaVaWa20}.  This uses a
generalization of bulge chasing called pole swapping and is a viable competitor
of QZ for the generalized eigenvalue problem.  It is natural to ask whether the
pole-swapping idea can be adapted to the standard eigenvalue problem for a
single matrix $A$, which would provide an alternative to Francis's algorithm.
The RQR algorithm introduced in this paper is such an alternative.

In Sections~\ref{sec:rqzrqr} and \ref{sec:rqr} we will explain the modifications
to RQZ that are needed to make an efficient and numerically stable RQR
algorithm, and in Section~\ref{sec:results} we will present some numerical
results that demonstrate that RQR is competitive with QR.  In fact, the RQR code
is slightly faster and a bit more accurate, in the sense that the backward error
is smaller.

\section{Modifying RQZ to make RQR}\label{sec:rqzrqr}

Ideally the reader should be familiar with some of the works about pole swapping cited above, especially the original
RQZ paper \cite{CaMeVa19a}.    However, even a reader who is unfamiliar 
with these developments will find this paper readable, as our presentation includes a high-level description
of RQZ.  For some fine points, recourse to the literature will be necessary.  

We will also use some core-chasing terminology, 
for which we refer to \cite{AuMaRoVaWa18} or \cite{AuMaVaWa15}, for example.   

The RQZ algorithm acts on \emph{Hessenberg pencils} $A - \lambda B$, that is,
pencils for which both $A$ and $B$ are upper-Hessenberg matrices.  Let's suppose
throughout that the matrices have dimension $n$.  The \emph{poles} of the
Hessenberg pencil $A - \lambda B$ are $\sigma_{1} = a_{21}/b_{21}$,
$\sigma_{2} = a_{32}/b_{32}$, $\dotsc$, $\sigma_{n-1} = a_{n,n-1}/b_{n,n-1}$,
i.e., the ratios of the subdiagonal elements.  Associated with this pencil is a
\emph{pole pencil} $A_{\pi} - \lambda B_{\pi}$ obtained by deleting the first
row and last column from $A - \lambda B$.  The pole pencil is clearly
$(n-1)\times (n-1)$ and upper triangular, and its eigenvalues are exactly the
poles of $A - \lambda B$.  Therefore the task of swapping two adjacent poles in
$A - \lambda B$ is exactly the same as that of interchanging two adjacent
eigenvalues of the pole pencil, a process that is well understood \cite{VanD81}.
The RQZ algorithm consists mainly of a large number of such swaps.

If we wish to find the eigenvalues of a single upper-Hessenberg matrix $A$, we
can do this by applying the RQZ algorithm to the pencil $A - \lambda I$, which
is a Hessenberg pencil.  However RQZ applies unitary equivalence transformations
that are not similarities, so the form $A - \lambda I$ is not preserved.  Once
we begin the iterations, the form of the pencil will become $A - \lambda U$,
where $U \neq I$.  However $U$ will be upper Hessenberg, and it will also be
unitary because $I$ is unitary and so are all of the transformations performed
by the RQZ algorithm.
Our task is just to show how to do the operations efficiently in this special
case.\footnote{We also investigated whether there is an efficient
      way of obtaining $A-\lambda U$ with $A$ upper Hessenberg and $U$ unitary
      upper Hessenberg from a general matrix $C$. However, we have so far not
      being able to find a better way than using the traditional reduction of
      $B$ to its Hessenberg form $A$ combined with $U=I$.}

It is well known, see \cite{AuMaRoVaWa18,AuMaVaWa15,m281,AuMaVaWa14b} for example, that a unitary Hessenberg matrix $U$
can be stored very compactly as a product of a descending sequence 
of $n-1$ core transformations:  $U = U_{1}U_{2} \cdots U_{n-1}$.
We remind the reader that a \emph{core transformation}  is just a unitary matrix whose active part is $2 \times 2$, 
acting on two successive rows and columns as indicated by its subscript.   For example, Givens rotations are core transformations.   Here $U_{1}$ acts
on rows 1 and 2, $U_{2}$ acts on rows 2 and 3, and so on.  We represent this pictorially by a double arrow pointing at the active rows. For example,
\begin{equation}\label{eq:uhess}
U \ \ = \ \ \parbox{2cm}{
\begin{tikzpicture}[scale=1.66,y=-1cm]
\foreach \j in {0,...,4}{\tikzrotation{-1.2+\j/5}{\j/5}}
\end{tikzpicture}
}.
\end{equation}
Here we have depicted the case $n=6$.  
Since each core
transformation is determined by two numbers, the total storage required for $U$ is just $2n-2$ numbers.   
  
The main operation in RQZ is the swap of two adjacent poles, which we call a move of type II.   
The interchange of poles $j-1$ and $j$ is accomplished by 
a transformation,
\begin{equation}
A \to Q_{j}^{*}AZ_{j-1}  \qquad U \to Q_{j}^{*}UZ_{j-1},
\label{eq:type2}
\end{equation}
where $Q_{j}$ and $Z_{j-1}$ are core transformations, $Q_{j}$ acting on rows $j$ and $j+1$ of $A$ and $U$, and $Z_{j-1}$ acting on columns $j-1$ and $j$.  We need to show how to carry out the transformation efficiently and accurately on $U$, 
which is stored as in \eqref{eq:uhess}.
There are several ways  to generate $Q_{j}$ and $Z_{j-1}$, and it's important how we do it here, but let's not worry about that
initially.   

In the case $j=3$ we have for the update of $U$ in \eqref{eq:type2} that,
\begin{displaymath}
Q_{j}^{*}UZ_{j-1} \ \ = \ \ \ Q_{3}^{*}UZ_{2} \ \ = \ \ \parbox{2cm}{
\begin{tikzpicture}[scale=1.66,y=-1cm]
\foreach \j in {0,...,4}{\tikzrotation{-1.2+\j/5}{\j/5}}
\tikzrotation[red]{-1.2}{.4}
\tikzrotation[blue]{-.6}{.2}x
\end{tikzpicture}
}.
\end{displaymath}
The result must again be upper Hessenberg, so there must be a way of reorganizing the core transformations into a  single descending
sequence, as pictured in \eqref{eq:uhess}.   One possibility is to move the blue core from the right to the left by a shift-through
operation (one turnover) \cite{AuMaRoVaWa18}, resulting in, 
\begin{displaymath}
Q_{j}^{*}UZ_{j-1} \ \ = \ \ \parbox{2.2cm}{
\begin{tikzpicture}[scale=1.66,y=-1cm]
\foreach \j in {0,...,4}{\tikzrotation{-1.2+\j/5}{\j/5}}
\tikzrotation[red]{-1.4}{.4}
\tikzrotation[gray]{-.6}{.2}
\tikzrotation[blue]{-1.2}{.4}
\turnoverrl{-.6}{.2}{-1.2}{.4};
\end{tikzpicture}
}.
\end{displaymath}
Because the matrix is upper Hessenberg, the resulting blue core must be the inverse of the red core, so that they cancel each other out.
Another possibility is to move the red core from the left to the right, resulting in,
\begin{displaymath}
Q_{j}^{*}UZ_{j-1} \ \ = \ \ \parbox{2cm}{
\begin{tikzpicture}[scale=1.66,y=-1cm]
\foreach \j in {0,...,4}{\tikzrotation{-1.2+\j/5}{\j/5}}
\tikzrotation[gray]{-1.2}{.4}
\turnoverlr{-1.2}{.4}{-.6}{.2}
\tikzrotation[red]{-.6}{.2}
\tikzrotation[blue]{-.4}{.2}
\end{tikzpicture}
}.
\end{displaymath}
Again, the resulting red core must be the inverse of the blue core, and they must cancel out.    Neither of these procedures
works in practice because of problems with roundoff errors, so we need to be a bit more careful.  

We have to go back to the details of the computation of $Q_{j}$ and $Z_{j-1}$,     
which must swap the eigenvalues of the subpencil, 
\begin{displaymath}
\left[\begin{array}{cc} a_{j,j-1} & a_{jj} \\ & a_{j+1,j} \end{array}\right] - \lambda
\left[\begin{array}{cc} u_{j,j-1} & u_{jj} \\ & u_{j+1,j} \end{array}\right].  
\end{displaymath}
For this we need to figure out the values of a few entries of $U$, which is,  fortunately, easy.  If the active part of the $j$th 
core transformation $U_j$ is $\left[\begin{smallmatrix} c_{j} & -s_{j} \\ s_{j} & \phantom{-}\overline{c}_{j} \end{smallmatrix}\right]$, 
then $u_{j,j-1} = s_{j-1}$ and 
$u_{jj} = \overline{c}_{j-1}c_{j}$, as the reader can easily check.   Thus our pencil is, 
\begin{equation}\label{eq:subpencil2}
\left[\begin{array}{cc} a_{j,j-1} & a_{jj} \\ & a_{j+1,j} \end{array}\right] - \lambda
\left[\begin{array}{cc} s_{j-1} & \overline{c}_{j-1}c_{j} \\ & s_{j} \end{array}\right]. 
\end{equation}
Its two eigenvalues are $\lambda_{1} = a_{j,j-1}/s_{j-1}$ 
and $\lambda_{2} = a_{j+1,j}/s_{j}$, and we need to swap them.   The basic procedure that we follow is due to 
Van Dooren \cite{VanD81}, as modified by Camps et al.\ \cite{CaMaVaWa20}.   There are several variants, and our 
choice of variant will depend on the relative magnitudes of the eigenvalues.  

\subsection*{The case $\absval{\lambda_{1}} \ge \absval{\lambda_{2}}$}

For this case we will use a variant that computes $Z_{j-1}$ first and then $Q_{j}$.  Substituting 
$\lambda_{2}$ for $\lambda$ in \eqref{eq:subpencil2}, we get,
\begin{equation*}
\begin{split}
H  & = s_{j}\left[\begin{array}{cc} a_{j,j-1} & a_{jj} \\ & a_{j+1,j} \end{array}\right] - a_{j+1,j}
\left[\begin{array}{cc} s_{j-1} & \overline{c}_{j-1}c_{j} \\ & s_{j} \end{array}\right], \\
& = \left[\begin{array}{cc}  s_{j}a_{j,j-1} -a_{j+1,j}s_{j-1} & s_{j}a_{jj} -
a_{j+1,j} \overline{c}_{j-1}c_{j}  \\ 0 & 0 \end{array}\right] = \left[\begin{array}{cc} {*} & {*} \\ 0 & 0 \end{array}\right].
\end{split}
\end{equation*}
Let $Z$ be a core transformation such that, 
\begin{displaymath}
HZ = \left[\begin{array}{cc} 0 & {*} \\ 0 & 0 \end{array}\right].
\end{displaymath}
This is (the active part of) our transformation $Z_{j-1}$ in \eqref{eq:type2}.   If we now compute $Q_{j}$ as prescribed in 
\cite{VanD81} or \cite{CaMaVaWa20}, we get a backward-stable swap that never fails, but since we are 
storing $U$ in the special form \eqref{eq:uhess}, we have to proceed differently.  

Instead of computing $Q_{j}$, we immediately apply
$Z_{j-1}$ to the pencil $A-\lambda U$.   In the case $j=3$, the picture looks like this:  
\begin{displaymath}
\parbox{3cm}{
\begin{tikzpicture}[scale=1.66,y=-1cm]
\draw (-.15,-.1) -- (-.2,-.1) -- (-.2,1.1) -- (-.15,1.1);
\draw (1.15,-.1) -- (1.2,-.1) -- (1.2,1.1) -- (1.15,1.1);
\foreach \j in {0,...,5}{
   \foreach \i in {\j,...,5}{\node at (\i/5,\j/5)
     [align=center,scale=1.0]{$\times$};}}
\foreach \j in {0,...,4}{\node at (\j/5,\j/5+.2)
     [align=center,scale=1.0]{$\times$};}
\tikzrotation[red]{1.4}{.2}
\end{tikzpicture}
} \qquad -\lambda \qquad
\parbox{2.2cm}{
\begin{tikzpicture}[scale=1.66,y=-1cm]
\foreach \j in {0,...,4}{\tikzrotation{\j/5}{\j/5}}
\tikzrotation[red]{.6}{.2};
\end{tikzpicture}
}.
\end{displaymath}
When we apply the core $Z_{j-1}$ to $A$ on the right, it recombines columns $2$ and $3$, creating a bulge in 
the $(4,2)$ position.   In $U$, let's do a turnover to move the extra core from the right to the left:
\begin{displaymath}
\parbox{3cm}{
\begin{tikzpicture}[scale=1.66,y=-1cm]
\draw (-.15,-.1) -- (-.2,-.1) -- (-.2,1.1) -- (-.15,1.1);
\draw (1.15,-.1) -- (1.2,-.1) -- (1.2,1.1) -- (1.15,1.1);
\foreach \j in {0,...,5}{
   \foreach \i in {\j,...,5}{\node at (\i/5,\j/5)
     [align=center,scale=1.0]{$\times$};}}
\foreach \j in {0,...,4}{\node at (\j/5,\j/5+.2)
     [align=center,scale=1.0]{$\times$};}
\node at (.2,.6)[align=center,scale=1.0,color=red]{$\times$};
\end{tikzpicture}
}  \qquad -\lambda \qquad
\parbox{2.2cm}{
\begin{tikzpicture}[scale=1.66,y=-1cm]
\foreach \j in {0,...,4}{\tikzrotation{\j/5}{\j/5}}
\tikzrotation[gray]{.6}{.2};
\tikzrotation[red]{0}{.4};
\turnoverrl{.6}{.2}{0}{.4};
\end{tikzpicture}
}.
\end{displaymath}
Now we need to apply $Q_{j}^{*}$ on the left to return the matrices to upper Hessenberg form.  How do we compute
$Q_{j}$?   The picture tells the story.   The remaining red core must be $Q_{j}$, as the way to return $U$ to Hessenberg form
is to get rid of the red core by multiplying by its inverse, which we show in blue here:  
\begin{displaymath}
\parbox{3cm}{
\begin{tikzpicture}[scale=1.66,y=-1cm]
\draw (-.15,-.1) -- (-.2,-.1) -- (-.2,1.1) -- (-.15,1.1);
\draw (1.15,-.1) -- (1.2,-.1) -- (1.2,1.1) -- (1.15,1.1);
\foreach \j in {0,...,5}{
   \foreach \i in {\j,...,5}{\node at (\i/5,\j/5)
     [align=center,scale=1.0]{$\times$};}}
\foreach \j in {0,...,4}{\node at (\j/5,\j/5+.2)
     [align=center,scale=1.0]{$\times$};}
\node at (.2,.6)[align=center,scale=1.0,color=red]{$\times$};
\tikzrotation[blue]{-.4}{.4}
\end{tikzpicture}
}  \qquad -\lambda \qquad
\parbox{2.2cm}{
\begin{tikzpicture}[scale=1.66,y=-1cm]
\foreach \j in {0,...,4}{\tikzrotation{\j/5}{\j/5}}
\tikzrotation[red]{0}{.4};
\tikzrotation[blue]{-.2}{.4}
\end{tikzpicture}
}.
\end{displaymath}
This returns $U$ to upper Hessenberg form, and it also knocks out the bulge in $A$, returning it to Hessenberg form:
\begin{displaymath}
\parbox{2.65cm}{
\begin{tikzpicture}[scale=1.66,y=-1cm]
\draw (-.15,-.1) -- (-.2,-.1) -- (-.2,1.1) -- (-.15,1.1);
\draw (1.15,-.1) -- (1.2,-.1) -- (1.2,1.1) -- (1.15,1.1);
\foreach \j in {0,...,5}{
   \foreach \i in {\j,...,5}{\node at (\i/5,\j/5)
     [align=center,scale=1.0]{$\times$};}}
\foreach \j in {0,...,4}{\node at (\j/5,\j/5+.2)
     [align=center,scale=1.0]{$\times$};}
\end{tikzpicture}
}  \qquad -\lambda\qquad
\parbox{2.2cm}{
\begin{tikzpicture}[scale=1.66,y=-1cm]
\foreach \j in {0,...,4}{\tikzrotation{\j/5}{\j/5}}
\end{tikzpicture}
}.
\end{displaymath}
This completes the swap.

This procedure is guaranteed to keep $U$ perfectly in Hessenberg form.  The entry in the $(j+1,j-1)$ position of $A$
will be slightly nonzero due to roundoff, but the error analysis in \cite{CaMaVaWa20}
guarantees that it is small enough to be ignored.   

\subsection*{The case $\absval{\lambda_{1}} < \absval{\lambda_{2}}$}

In this case we compute $Q_{j}$ first.  Substituting $\lambda_{1}$ for $\lambda$ in \eqref{eq:subpencil2},
we get, 
\begin{equation*}
\begin{split}
H  & = s_{j-1}\left[\begin{array}{cc} a_{j,j-1} & a_{jj} \\ & a_{j+1,j} \end{array}\right] - a_{j,j-1}
\left[\begin{array}{cc} s_{j-1} & \overline{c}_{j-1}c_{j} \\ & s_{j} \end{array}\right], \\
& = \left[\begin{array}{cc}  0 & s_{j-1}a_{jj} - a_{j,j-1}\overline{c}_{j-1}c_{j} \\ 0 & s_{j-1}a_{j+1,j} -a_{j,j-1}s_{j}\end{array}\right] 
= \left[\begin{array}{cc} 0 & {*} \\ 0 & {*} \end{array}\right].
\end{split}
\end{equation*}
Now compute a core transformation $Q$ such that, 
\begin{displaymath}
Q^{*}H = \left[\begin{array}{cc}  0 & {*} \\ 0 & 0 \end{array}\right].
\end{displaymath}
This $Q$ is (the active part of) our 
desired $Q_{j}$ in \eqref{eq:type2}.  Now we apply $Q_{j}^{*}$ to the pencil
$A-\lambda U$ immediately.  
In the case $j=3$ it looks like this:
\begin{displaymath}
\parbox{3cm}{
\begin{tikzpicture}[scale=1.66,y=-1cm]
\draw (-.15,-.1) -- (-.2,-.1) -- (-.2,1.1) -- (-.15,1.1);
\draw (1.15,-.1) -- (1.2,-.1) -- (1.2,1.1) -- (1.15,1.1);
\foreach \j in {0,...,5}{
   \foreach \i in {\j,...,5}{\node at (\i/5,\j/5)
     [align=center,scale=1.0]{$\times$};}}
\foreach \j in {0,...,4}{\node at (\j/5,\j/5+.2)
     [align=center,scale=1.0]{$\times$};}
\tikzrotation[red]{-.4}{.4}
\end{tikzpicture}
}  \qquad   -\lambda\qquad
\parbox{2.2cm}{
\begin{tikzpicture}[scale=1.66,y=-1cm]
\foreach \j in {0,...,4}{\tikzrotation{\j/5}{\j/5}}
\tikzrotation[red]{0}{.4}
\end{tikzpicture}
}.
\end{displaymath}
 When we apply $Q_{j}^{*}$ to $A$, it recombines rows 3 and 4, making a bulge in the $(4,2)$ position.   We 
 pass $Q_{j}^{*}$ through $U$ by a turnover:
 \begin{displaymath}
\parbox{3cm}{
\begin{tikzpicture}[scale=1.66,y=-1cm]
\draw (-.15,-.1) -- (-.2,-.1) -- (-.2,1.1) -- (-.15,1.1);
\draw (1.15,-.1) -- (1.2,-.1) -- (1.2,1.1) -- (1.15,1.1);
\foreach \j in {0,...,5}{
   \foreach \i in {\j,...,5}{\node at (\i/5,\j/5)
     [align=center,scale=1.0]{$\times$};}}
\foreach \j in {0,...,4}{\node at (\j/5,\j/5+.2)
     [align=center,scale=1.0]{$\times$};}
\node at (.2,.6)[align=center,scale=1.0,color=red]{$\times$};
\end{tikzpicture}
}  \qquad  -\lambda \qquad
\parbox{2.2cm}{
\begin{tikzpicture}[scale=1.66,y=-1cm]
\foreach \j in {0,...,4}{\tikzrotation{\j/5}{\j/5}}
\tikzrotation[red]{.6}{.2};
\tikzrotation[gray]{0}{.4};
\turnoverlr{0}{.4}{.6}{.2};
\end{tikzpicture}
}.
\end{displaymath}
We now return $U$ to upper Hessenberg form by multiplying by the inverse of the extra core, which we 
mark in blue.  This must be $Z_{j-1}$.
\begin{displaymath}
\parbox{3cm}{
\begin{tikzpicture}[scale=1.66,y=-1cm]
\draw (-.15,-.1) -- (-.2,-.1) -- (-.2,1.1) -- (-.15,1.1);
\draw (1.15,-.1) -- (1.2,-.1) -- (1.2,1.1) -- (1.15,1.1);
\foreach \j in {0,...,5}{
   \foreach \i in {\j,...,5}{\node at (\i/5,\j/5)
     [align=center,scale=1.0]{$\times$};}}
\foreach \j in {0,...,4}{\node at (\j/5,\j/5+.2)
     [align=center,scale=1.0]{$\times$};}
\node at (.2,.6)[align=center,scale=1.0,color=red]{$\times$};
\tikzrotation[blue]{1.4}{.2}
\end{tikzpicture}
} \qquad -\lambda\qquad
\parbox{2.2cm}{
\begin{tikzpicture}[scale=1.66,y=-1cm]
\foreach \j in {0,...,4}{\tikzrotation{\j/5}{\j/5}}
\tikzrotation[red]{.6}{.2};
\tikzrotation[blue]{.8}{.2};
\end{tikzpicture}
}.
\end{displaymath}
When we apply $Z_{j-1}$ to $A$ on the right, it recombines columns 2 and 3, cancelling out the bulge.
\begin{displaymath}
\parbox{3cm}{
\begin{tikzpicture}[scale=1.66,y=-1cm]
\draw (-.15,-.1) -- (-.2,-.1) -- (-.2,1.1) -- (-.15,1.1);
\draw (1.15,-.1) -- (1.2,-.1) -- (1.2,1.1) -- (1.15,1.1);
\foreach \j in {0,...,5}{
   \foreach \i in {\j,...,5}{\node at (\i/5,\j/5)
     [align=center,scale=1.0]{$\times$};}}
\foreach \j in {0,...,4}{\node at (\j/5,\j/5+.2)
     [align=center,scale=1.0]{$\times$};}
\end{tikzpicture}
}  \qquad-\lambda \qquad
\parbox{2.2cm}{
\begin{tikzpicture}[scale=1.66,y=-1cm]
\foreach \j in {0,...,4}{\tikzrotation{\j/5}{\j/5}}
\end{tikzpicture}
}
\end{displaymath}
This completes the swap.

Again the $(j+1,j-1)$ entry of $A$ will be slightly nonzero due to roundoff, but it is guaranteed to be 
small enough to ignore.  The analysis in \cite{CaMaVaWa20} does not mention this case explicitly, but 
this is dual to the ``$Z$-first'' method shown above and has the same numerical properties.  

\subsection*{Moves of type I}

Each iteration of the RQR algorithm begins and ends with a move of type I.   
A type I move at the top of the pencil inserts an arbitrary pole $\rho$, replacing $\sigma_{1}$, at the 
top of the pencil.   Let $x = (A - \rho B)e_{1}$, and let $Q_{1}$ be a core transformation such that 
$Q_{1}^{*}x = \alpha e_{1}$ for some $\alpha \neq 0$.   This is possible because only the first two entries of $x$
are nonzero:  $\left[\begin{smallmatrix} x_{1} \\ x_{2} \end{smallmatrix}\right] = 
\left[\begin{smallmatrix} a_{11} - \rho u_{11} 
\\ a_{21} - \rho u_{21} \end{smallmatrix}\right]$.   Since $U$ is stored in the form \eqref{eq:uhess}, we need to extract
the values of $u_{11}$ and $u_{21}$, but this is easy:   if 
$\left[\begin{smallmatrix} c_{1} & -s_{1} \\ s_{1} & \phantom{-}\overline{c}_{1} \end{smallmatrix}\right]$ is the active part
of the first core transformation in $U$, then $u_{11} = c_{1}$ and $u_{21} = s_{1}$.   

Once we have $Q_{1}^{*}$, we apply it to $A-\lambda U$ to obtain $Q_{1}^{*}A -
\lambda Q_{1}^{*}U$:
\begin{displaymath}
\parbox{3cm}{
\begin{tikzpicture}[scale=1.66,y=-1cm]
\draw (-.15,-.1) -- (-.2,-.1) -- (-.2,1.1) -- (-.15,1.1);
\draw (1.15,-.1) -- (1.2,-.1) -- (1.2,1.1) -- (1.15,1.1);
\foreach \j in {0,...,5}{
   \foreach \i in {\j,...,5}{\node at (\i/5,\j/5)
     [align=center,scale=1.0]{$\times$};}}
\foreach \j in {0,...,4}{\node at (\j/5,\j/5+.2)
     [align=center,scale=1.0]{$\times$};}
\tikzrotation[red]{-.4}{0}
\end{tikzpicture}
}  \qquad -\lambda\qquad
\parbox{2.2cm}{
\begin{tikzpicture}[scale=1.66,y=-1cm]
\foreach \j in {0,...,4}{\tikzrotation{\j/5}{\j/5}}
\tikzrotation[red]{-.2}{0}
\end{tikzpicture}
}.
\end{displaymath}
The core $Q_{1}^{*}$ recombines the first two rows of $A$, and  it can be absorbed into
$U$ by a fusion operation.   The resulting $Q_{1}^{*}A -\lambda Q_{1}^{*}U$ remains in Hessenberg form.
This completes the move of type I.   The reader can easily check that the first pole is
now $\rho$.   The other poles are unchanged.   

A move of type I at the bottom replaces the bottom pole by an arbitrary pole $\tau$.   
The details of this routine are analogous and are left to the reader.

\section{The RQR algorithm}\label{sec:rqr}

The RQR algorithm applied to  $A-\lambda U$ is the same as the RQZ algorithm, with the modifications
described in the previous section.   Each iteration of the basic algorithm begins with the choice of 
a shift $\rho$.   Any shifting strategy that is commonly used by the Francis or Moler-Stewart algorithm can be
used.  For example, the simplest choice is the \emph{Rayleigh-quotient shift}  $\rho = a_{nn}/u_{nn}$.  
A better choice is the \emph{Wilkinson shift}, which computes the
eigenvalues of the $2 \times 2$ subpencil in the lower right-hand corner and
takes $\rho$ to be the eigenvalue that is closer to $a_{nn}/u_{nn}$.

Once $\rho$ has been chosen, it is inserted as a pole at the top of the pencil by a move of type I.  
Then, by a sequence of moves of type II, $\rho$ is swapped with poles $\sigma_{2}$, $\sigma_{3}$, and 
so on, until $\rho$ reaches the bottom of the pencil.   Then it is replaced by a new pole $\tau$ by a move
of type I.   The simplest pole choice is the \emph{Rayleigh-quotient pole} $\tau = a_{11}/u_{11}$.   A better
choice is the \emph{Wilkinson pole}:   Compute the eigenvalues of the $2 \times 2$ submatrix in the upper
left-hand corner and take $\tau$ to be the eigenvalue that is closer to $a_{11}/u_{11}$.  
This completes the iteration.   

Repeated iterations will generally cause the pencil to converge to triangular form, revealing the eigenvalues
on the main diagonal.   Convergence is not equally rapid up and down the pencil.   Good shifts $\rho$ inserted
at the top will cause rapid convergence at the bottom.    If $a_{n,n-1} \to 0$ and $u_{n,n-1} \to 0$ rapidly, they 
can be declared to be zero after just a few iterations, and $a_{nn}/u_{nn}$ can be declared to be an eigenvalue.  
This is known as a deflation and it reduces the problem to size $n-1$. After deflating, we can go after the next eigenvalue.
At the same time, the poles 
$\tau$ that are inserted at the bottom will gradually move to the top and improve the rate of convergence to zero 
of elements at the top, such as $a_{21}$ and $b_{21}$.   
This generally decreases the total number of iterations required compared to a similar bulge chasing approach. 

The algorithm can be enhanced in several ways.  For one thing, there is no need
to wait until one iteration is complete before initiating the next one; many
iterations can proceed at once.  Suppose we want to do $m$ iterations
simultaneously, where $1 \ll m \ll n$.  We pick $m$ shift
$\rho_{1},\ldots,\rho_{m}$ at once.  For example, we can compute the eigenvalues
of the $m\times m$ subpencil in the lower right-hand corner of $A -\lambda B$
and use them as shifts.  Then, by moves of types I and II, we insert
$\rho_{1},\ldots,\rho_{m}$ as the first $m$ poles in the pencil.  Then we chase
them all down to the bottom simultaneously.  This improves performance
substantially by enabling the use of level-3 BLAS and decreasing cache misses.
This was proposed in the context of bulge chasing by Braman et al.\
\cite{BrByMa01} and Lang \cite{Lang98}, and discussed in the context of pole
swapping in \cite{CaMaVaWa20,StCaMeVa21}.

Another very important enhancement is the use of aggressive early deflation,
introduced by Braman et al.\ \cite{BrByMa01b} and discussed in the context of
pole swapping in \cite{StCaMeVa21}.

For real matrices it is desirable to introduce a double-shift algorithm that allows for the use of complex-conjugate 
shifts and poles while staying in real arithmetic.   

We have not implemented any of these enhancements so far.  

\section{The Difference between QR and RQR}

The QZ algorithm acts on pencils $A - \lambda B$ in Hessenberg-triangular form:  $A$ is upper Hessenberg, and 
$B$ is upper triangular.    This is a special type of Hessenberg pencil in which all of the poles are infinite.  
In \cite{CaMaVaWa20}  we showed that if the RQZ algorithm is applied to a Hessenberg-triangular pencil, it reduces
to the QZ algorithm.   This equivalence depends on inserting an infinite pole at the bottom (the final move of Type I at the
end of the iteration), which is necessary to maintain the upper-triangular form of $B$.  

The RQR algorithm developed here begins with a pencil $A - \lambda I$, which is a Hessenberg-triangular pencil.
The reader can check that if the QZ algorithm is applied to such a pencil it reduces essentially to QR.  The pencil 
may change to $A - \lambda U$, but $U$ is always unitary and diagonal, which means that $U$ is nearly the identity 
matrix.   With appropriate normalizations the form $U=I$ can be preserved, in which case QZ really does reduce to
QR.     Thus it may appear that RQZ, which reduces to QZ, is really just QR in this special case.    
So what is the difference?   

The difference is that at the end of an RQR step we are not obliged to insert an infinite pole, because we are not
interested in preserving the diagonal form of U.   Instead we insert a finite pole chosen with an eye
to improving convergence at the top of the pencil.   This is the only difference, and this difference is significant, as the following results show.

\section{Numerical Results}\label{sec:results}

We wrote a Fortran implementation of the RQR algorithm with Wilkinson shifts and Wilkinson poles.   
We compared the performance of our RQR code (ZLAHPS) with that of ZLAHQR, the QR
kernel from LAPACK version 3.12.0.  ZLAHPS was built by modifying
ZLAHQR.\footnote{We also modified 
routines for manipulating rotations, for instance the turnover, from the
eiscor project \cite{eiscor}.}  On
the test problems that we considered, RQR was, on average, faster and more
accurate than QR.

For the numerical experiments reported below we used the gfortran compiler from gcc 9.4.0 with
optimization flag -O2 on a computer with Ubuntu 20.04.2, an Intel Core i7-10700K
CPU with 12 MiB of L3 cache, and 32 GiB of RAM. Experiments on other computers
showed similar results.

\begin{figure}[htp]
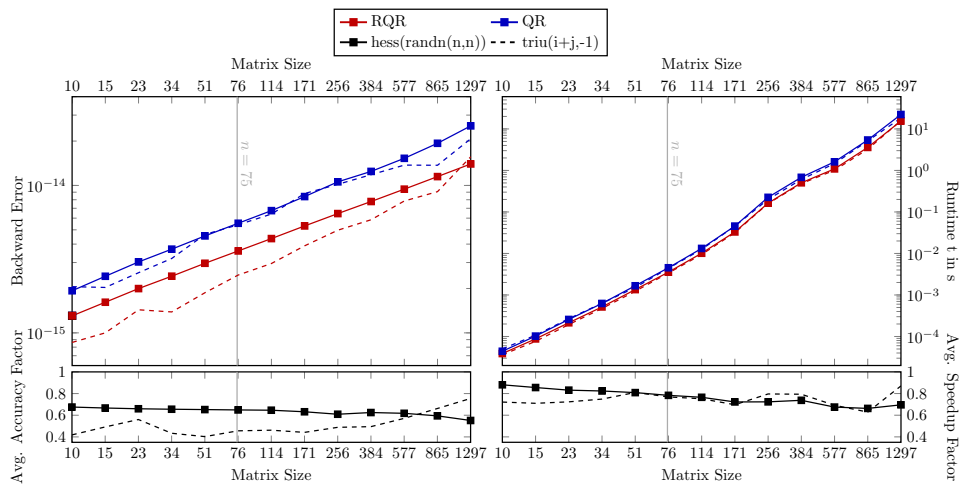

  \begin{tikzpicture}[scale=0.6]
  \begin{semilogyaxis}[%
    name = bwe1,
    width=0.80\textwidth,
    height=0.36\textheight,
    xticklabel pos = lower,%
    xticklabel pos = upper,%
    ylabel = {Backward Error},%
    xlabel = {Matrix Size},%
    xmin = 1,%
    xmax = 13,%
    ymin = 6e-16,%
    ymax = 4e-14,%
    every axis y label/.style= {at={(ticklabel cs:0.5)},rotate=90,yshift=3pt,
    },%
    xtick={1,2,3,4,5,6,7,8,9,10,11,12,13},%
    xticklabels={10, 15, 23, 34, 51, 76, 114, 171, 256, 384, 577, 865, 1297},%
    every axis legend/.append style={ at={(1.0,1.15)}, anchor=south,
      cells={anchor=west},
      legend columns=2,
      thick},%
    thick
    ]
    
    \draw[SPECgray!50] (axis cs:5.9697,6e-16)--node [near end,rotate = -90,yshift=6pt] {$n=75$}(axis cs:5.9697,4e-14);

    \addplot[%
    SPECred,%
    mark=square*,%
    mark options = {solid, fill = SPECred},%
    very thick,%
    ]coordinates{%
      (1,1.304240e-15)%
    };
    \addplot[%
    SPECblue,%
    mark=square*,%
    mark options = {solid, fill = SPECblue},%
    very thick,%
    ]coordinates{%
      (1,1.304240e-15)%
    };
    \addplot[%
    black,%
    mark=square*,%
    mark options = {solid, fill = black},%
    very thick,%
    ]coordinates{%
      (1,1.304240e-15)%
    };
    \addplot[%
    black,%
    dashed,%
    mark options = {solid, fill = black},%
    very thick,%
    ]coordinates{%
      (1,1.304240e-15)%
    };
    
    \input{bwe.tex}
    
    \legend{RQR,QR,{hess(randn(n,n))},{triu(i+j,-1)}}
    
  \end{semilogyaxis}
  \begin{axis}[
    name = bwe2,
    at = {($(bwe1.south)+(0,-4pt)$)},
    anchor = north,
    width=0.80\textwidth,
    height=0.15\textheight,
    yticklabel pos = left,%
    xticklabel pos = lower,%
    ymin = 0.35,
    ymax = 1,
    xmin = 1,
    xmax = 13,
    ylabel = {Avg. Accuracy Factor},%
    xlabel = {Matrix Size},%
    every axis y label/.style= {at={(ticklabel cs:0.5)},rotate=90,yshift=14pt,
    },%
    xtick={1,2,3,4,5,6,7,8,9,10,11,12,13},%
    xticklabels={10, 15, 23, 34, 51, 76, 114, 171, 256, 384, 577, 865, 1297},%
    every axis legend/.append style={ at={(0.5,1.05)}, anchor=south,
      cells={anchor=west},
      legend columns=3},%
    thick
    ]

    \draw[SPECgray!50] (axis cs:5.9697,0.35)--(axis cs:5.9697,1);

    \input{bwe_ratio.tex}
    
  \end{axis}
  
  \begin{semilogyaxis}[
    name = runtime1,
    at = {($(bwe1.east)+(+20pt,0)$)},
    anchor = west,
    width=0.80\textwidth,
    height=0.36\textheight,
    xticklabel pos = upper,%
    yticklabel pos = right,%
    ylabel = {Runtime t in s},%
    xlabel = {Matrix Size},%
    xmin = 1,%
    xmax = 13,%
    ymin = 2e-5,%
    ymax = 60,%
    every axis y label/.style= {at={(ticklabel cs:0.5)},rotate=-90,yshift=5pt,
    },%
    xtick={1,2,3,4,5,6,7,8,9,10,11,12,13},%
    xticklabels={10, 15, 23, 34, 51, 76, 114, 171, 256, 384, 577, 865, 1297},%
    every axis legend/.append style={ at={(0.5,1.05)}, anchor=south,
      cells={anchor=west},
      legend columns=2},%
    thick
    ]
    
    \draw[SPECgray!50] (axis cs:5.9697,2e-5)--node [near end,rotate = -90,yshift=6pt] {$n=75$}(axis cs:5.9697,60);

    \input{runtime.tex}
    
    
  \end{semilogyaxis}
  \begin{axis}[
    name = runtime2,
    at = {($(runtime1.south)+(0,-4pt)$)},
    anchor = north,
    width=0.80\textwidth,
    height=0.15\textheight,
    yticklabel pos = right,%
    xticklabel pos = lower,%
    ymin = 0.35,
    ymax = 1,
    xmin = 1,
    xmax = 13,
    ylabel = {Avg. Speedup Factor},%
    xlabel = {Matrix Size},%
    every axis y label/.style= {at={(ticklabel cs:0.5)},rotate=-90,yshift=12pt,
    },%
    xtick={1,2,3,4,5,6,7,8,9,10,11,12,13},%
    xticklabels={10, 15, 23, 34, 51, 76, 114, 171, 256, 384, 577, 865, 1297},%
    every axis legend/.append style={ at={(0.5,1.05)}, anchor=south,
      cells={anchor=west},
      legend columns=3},%
    thick
    ]

    \draw[SPECgray!50] (axis cs:5.9697,0.35)--(axis cs:5.9697,1);
    
    \input{runtime_ratio.tex}
    
  \end{axis}
\end{tikzpicture}

  \caption{Comparison of the RQR pole swapping algorithm (in
    \textcolor{SPECred}{red}) with the QR bulge chasing algorithm (in
    \textcolor{SPECblue}{blue}) in terms of backward error (left) and runtime
    (right) for random matrices reduced to upper Hessenberg form (square
    markers) and $i + j$ upper Hessenberg matrices (dashed line).}
\label{fig}
\end{figure}

We considered two kinds of matrices: (1) randomly-generated matrices reduced to
upper Hessenberg form,  and (2) upper Hessenberg matrices with entries
$a_{ij} = i + j$ ($i < j+2$).  Matrices of dimension $10 \times 10$ up to
$1297\times 1297$ were tested.  For the random matrices we did 100 trials at
each dimension.  The results are shown in Figure~\ref{fig}.

The RQR algorithm (red) is consistently faster than QR (blue).
The typical speedup is around 17\% for matrices of size less than $75$ and 29\%
for larger matrices. At the same time, the backward errors are smaller by a
factor of 1.5 to~2. The results for smaller matrices are particularly relevant
since LAPACK's main eigenvalue routine ZHSEQR uses ZLAHQR only for matrices of dimension up
to~75.\footnote{The LAPACK installation used for the experiment uses
  $n_{\min}=75$ for deciding when to use ZLAHQR. This parameter is machine- and
  installation-dependent, thus your mileage may vary.}

The results for the randomly-generated matrices are also given in tabular form in Table~\ref{tab1}.
This lists execution times, backward errors (BWE), and total iterations divided by $n$ (It/n).
The iteration counts are significantly less for RQR, which explains why RQR is 
faster, and may also explain why it is more accurate, as fewer iterations imply fewer
roundoff errors.  


We also tested both codes on 35 matrices of size less than 1000
from matrix market \cite{matrixmarket}.   Both codes performed well.
With a  small number of exceptions, RQR was faster and more accurate
than QR.  See Table~\ref{tab2}.

These results show that pole swapping is a viable alternative to bulge chasing
for the standard eigenvalue problem. The code and the experiments we presented
only use a single shift in each iteration. This is the natural and the
fastest choice for small matrices of dimension up to about 75.
For larger matrices multishift/multibulge variants with aggressive early
deflation are preferred. We do not yet have pole swapping code with these
features. Hence, a fair comparison with LAPACK's routines including these
features is not possible at this time. Nevertheless, the experiments that we 
have presented here 
are a preliminary indication that pole swapping has the potential to be faster
and more accurate than the QR codes that are currently in use. 

\bibliographystyle{siam}
\bibliography{rqr}

\begin{table}[tbp]
  \centering
  \begin{tabular}{rccc|rccc}
\toprule
\multicolumn{3}{l}{RQR}&&\multicolumn{3}{l}{QR}&\\
\midrule
n & Time [s] & BWE & It/n &
n & Time [s] & BWE & It/n\\
\midrule
10 %
&   3.88e-05 
&   1.30e-15 
&   2.58 
& 10 %
&   4.41e-05 
&   1.93e-15 
&   3.10 
\\
15 %
&   8.74e-05 
&   1.61e-15 
&   2.67 
& 15 %
&   1.02e-04 
&   2.42e-15 
&   3.27 
\\
23 %
&   2.14e-04 
&   2.00e-15 
&   2.70 
& 23 %
&   2.57e-04 
&   3.03e-15 
&   3.35 
\\
34 %
&   5.15e-04 
&   2.42e-15 
&   2.74 
& 34 %
&   6.25e-04 
&   3.70e-15 
&   3.38 
\\
51 %
&   1.35e-03 
&   2.97e-15 
&   2.74 
& 51 %
&   1.67e-03 
&   4.55e-15 
&   3.37 
\\
76 %
&   3.54e-03 
&   3.59e-15 
&   2.74 
& 76 %
&   4.52e-03 
&   5.54e-15 
&   3.35 
\\
114 %
&   1.01e-02 
&   4.36e-15 
&   2.74 
& 114 %
&   1.32e-02 
&   6.75e-15 
&   3.32 
\\
171 %
&   3.30e-02 
&   5.31e-15 
&   2.73 
& 171 %
&   4.57e-02 
&   8.41e-15 
&   3.28 
\\
256 %
&   1.63e-01 
&   6.45e-15 
&   2.72 
& 256 %
&   2.25e-01 
&   1.06e-14 
&   3.25 
\\
384 %
&   5.03e-01 
&   7.78e-15 
&   2.70 
& 384 %
&   6.83e-01 
&   1.25e-14 
&   3.22 
\\
577 %
&   1.09 
&   9.44e-15 
&   2.68 
& 577 %
&   1.61 
&   1.53e-14 
&   3.19 
\\
865 %
&   3.59 
&   1.15e-14 
&   2.67 
& 865 %
&   5.41 
&   1.93e-14 
&   3.16 
\\
1297 %
&  15.51 
&   1.40e-14 
&   2.66 
& 1297 %
&  22.30 
&   2.54e-14 
&   3.13 
\\
    \bottomrule
  \end{tabular}
  \caption{Comparison of RQR vs.\ QR on random matrices with
    normally-distributed entries. In a preprocessing step the matrices have been
    reduced to upper Hessenberg form.}
  \label{tab1}
\end{table}

\begin{table}[p]
\begin{center}
  \begin{tabular}{ll|ccc|ccc}
\toprule
&&\multicolumn{3}{l|}{RQR}&\multicolumn{3}{l}{QR}\\
    n & Name & Time [s] & BWE 
    & It/n
    & Time [s] & BWE 
                            & It/n\\
\midrule
200 &bwm200 &{\color{SPECgreen!60!black} 0.031} &{\color{SPECred}1.19e-14} &
{\color{SPECgreen!60!black} 1.81} 
    &{\color{SPECred} 0.037} &{\color{SPECgreen!60!black}6.27e-15} &
{\color{SPECred} 2.04} 
\\
961 &cdde1 &{\color{SPECgreen!60!black} 3.460} &{\color{SPECgreen!60!black}8.17e-15} &
{\color{SPECgreen!60!black} 1.89} 
    &{\color{SPECred} 3.779} &{\color{SPECred}1.18e-14} &
{\color{SPECred} 2.02} 
\\
961 &cdde2 &{\color{SPECgreen!60!black} 4.791} &{\color{SPECgreen!60!black}9.88e-15} &
{\color{SPECgreen!60!black} 2.53} 
    &{\color{SPECred} 6.706} &{\color{SPECred}1.61e-14} &
{\color{SPECred} 3.07} 
\\
961 &cdde3 &{\color{SPECgreen!60!black} 3.345} &{\color{SPECgreen!60!black}7.71e-15} &
{\color{SPECgreen!60!black} 1.90} 
    &{\color{SPECred} 3.787} &{\color{SPECred}1.18e-14} &
{\color{SPECred} 2.03} 
\\
961 &cdde4 &{\color{SPECgreen!60!black} 4.962} &{\color{SPECgreen!60!black}1.00e-14} &
{\color{SPECgreen!60!black} 2.53} 
    &{\color{SPECred} 6.879} &{\color{SPECred}1.61e-14} &
{\color{SPECred} 3.08} 
\\
961 &cdde5 &{\color{SPECgreen!60!black} 3.335} &{\color{SPECred}1.31e-14} &
{\color{SPECgreen!60!black} 1.91} 
    &{\color{SPECred} 3.769} &{\color{SPECgreen!60!black}1.25e-14} &
{\color{SPECred} 2.01} 
\\
961 &cdde6 &{\color{SPECgreen!60!black} 4.909} &{\color{SPECgreen!60!black}1.01e-14} &
{\color{SPECgreen!60!black} 2.52} 
    &{\color{SPECred} 6.753} &{\color{SPECred}1.58e-14} &
{\color{SPECred} 3.09} 
\\
104 &ck104 &{\color{SPECgreen!60!black} 0.005} &{\color{SPECgreen!60!black}2.25e-15} &
{\color{SPECred} 2.06} 
    &{\color{SPECred} 0.005} &{\color{SPECred}3.80e-15} &
{\color{SPECgreen!60!black} 2.03} 
\\
400 &ck400 &{\color{SPECgreen!60!black} 0.233} &{\color{SPECred}8.61e-15} &
{\color{SPECgreen!60!black} 1.82} 
    &{\color{SPECred} 0.320} &{\color{SPECgreen!60!black}7.91e-15} &
{\color{SPECred} 1.89} 
\\
656 &ck656 &{\color{SPECgreen!60!black} 1.037} &{\color{SPECgreen!60!black}9.14e-15} &
{\color{SPECred} 2.81} 
    &{\color{SPECred} 1.339} &{\color{SPECred}1.06e-14} &
{\color{SPECgreen!60!black} 2.00} 
\\
512 &dwa512 &{\color{SPECgreen!60!black} 0.647} &{\color{SPECgreen!60!black}1.05e-14} &
{\color{SPECgreen!60!black} 1.89} 
    &{\color{SPECred} 0.795} &{\color{SPECred}1.21e-14} &
{\color{SPECred} 1.91} 
\\
512 &dwb512 &{\color{SPECgreen!60!black} 0.618} &{\color{SPECgreen!60!black}5.06e-15} &
{\color{SPECgreen!60!black} 1.88} 
    &{\color{SPECred} 0.779} &{\color{SPECred}7.00e-15} &
{\color{SPECred} 1.93} 
\\
163 &lop163 &{\color{SPECgreen!60!black} 0.024} &{\color{SPECgreen!60!black}4.10e-15} &
{\color{SPECgreen!60!black} 2.37} 
    &{\color{SPECred} 0.029} &{\color{SPECred}7.83e-15} &
{\color{SPECred} 2.59} 
\\
100 &olm100 &{\color{SPECgreen!60!black} 0.003} &{\color{SPECgreen!60!black}2.72e-15} &
{\color{SPECgreen!60!black} 1.70} 
    &{\color{SPECred} 0.004} &{\color{SPECred}4.62e-15} &
{\color{SPECred} 2.02} 
\\
500 &olm500 &{\color{SPECgreen!60!black} 0.243} &{\color{SPECred}8.21e-15} &
{\color{SPECgreen!60!black} 1.63} 
    &{\color{SPECred} 0.262} &{\color{SPECgreen!60!black}7.53e-15} &
{\color{SPECred} 1.85} 
\\
225 &pde225 &{\color{SPECgreen!60!black} 0.068} &{\color{SPECgreen!60!black}4.83e-15} &
{\color{SPECgreen!60!black} 2.77} 
    &{\color{SPECred} 0.086} &{\color{SPECred}7.70e-15} &
{\color{SPECred} 3.04} 
\\
900 &pde900 &{\color{SPECgreen!60!black} 3.575} &{\color{SPECgreen!60!black}9.62e-15} &
{\color{SPECgreen!60!black} 2.49} 
    &{\color{SPECred} 4.900} &{\color{SPECred}1.48e-14} &
{\color{SPECred} 2.70} 
\\
362 &plat362 &{\color{SPECgreen!60!black} 0.202} &{\color{SPECgreen!60!black}6.93e-15} &
{\color{SPECgreen!60!black} 2.14} 
    &{\color{SPECred} 0.235} &{\color{SPECred}1.20e-14} &
{\color{SPECred} 2.38} 
\\
362 &plskz362 &{\color{SPECgreen!60!black} 0.255} &{\color{SPECgreen!60!black}7.42e-15} &
{\color{SPECgreen!60!black} 2.75} 
    &{\color{SPECred} 0.333} &{\color{SPECred}1.22e-14} &
{\color{SPECred} 3.25} 
\\
324 &qc324 &{\color{SPECgreen!60!black} 0.163} &{\color{SPECgreen!60!black}5.77e-15} &
{\color{SPECgreen!60!black} 2.22} 
    &{\color{SPECred} 0.220} &{\color{SPECred}9.72e-15} &
{\color{SPECred} 2.39} 
\\
768 &qh768 &{\color{SPECgreen!60!black} 1.216} &{\color{SPECred}3.63e-15} &
{\color{SPECred} 2.93} 
    &{\color{SPECred} 1.427} &{\color{SPECgreen!60!black}3.04e-15} &
{\color{SPECgreen!60!black} 2.11} 
\\
882 &qh882 &{\color{SPECgreen!60!black} 1.127} &{\color{SPECgreen!60!black}4.68e-15} &
{\color{SPECred} 3.70} 
    &{\color{SPECred} 1.214} &{\color{SPECred}5.56e-15} &
{\color{SPECgreen!60!black} 2.08} 
\\
480 &rbs480a &{\color{SPECgreen!60!black} 0.666} &{\color{SPECgreen!60!black}8.73e-15} &
{\color{SPECgreen!60!black} 2.62} 
    &{\color{SPECred} 0.853} &{\color{SPECred}1.35e-14} &
{\color{SPECred} 3.01} 
\\
480 &rbs480b &{\color{SPECgreen!60!black} 0.656} &{\color{SPECgreen!60!black}8.81e-15} &
{\color{SPECgreen!60!black} 2.70} 
    &{\color{SPECred} 1.066} &{\color{SPECred}1.58e-14} &
{\color{SPECred} 3.06} 
\\
200 &rdb200 &{\color{SPECgreen!60!black} 0.025} &{\color{SPECgreen!60!black}3.85e-15} &
{\color{SPECred} 2.21} 
    &{\color{SPECred} 0.028} &{\color{SPECred}6.26e-15} &
{\color{SPECgreen!60!black} 1.68} 
\\
200 &rdb200l &{\color{SPECgreen!60!black} 0.025} &{\color{SPECgreen!60!black}4.02e-15} &
{\color{SPECred} 2.78} 
    &{\color{SPECred} 0.026} &{\color{SPECred}5.66e-15} &
{\color{SPECgreen!60!black} 1.74} 
\\
450 &rdb450 &{\color{SPECgreen!60!black} 0.270} &{\color{SPECgreen!60!black}9.19e-15} &
{\color{SPECred} 2.40} 
    &{\color{SPECred} 0.309} &{\color{SPECred}1.14e-14} &
{\color{SPECgreen!60!black} 1.61} 
\\
450 &rdb450l &{\color{SPECgreen!60!black} 0.273} &{\color{SPECgreen!60!black}6.44e-15} &
{\color{SPECred} 2.52} 
    &{\color{SPECred} 0.303} &{\color{SPECred}9.64e-15} &
{\color{SPECgreen!60!black} 1.66} 
\\
800 &rdb800l &{\color{SPECgreen!60!black} 1.904} &{\color{SPECgreen!60!black}9.69e-15} &
{\color{SPECred} 2.38} 
    &{\color{SPECred} 2.164} &{\color{SPECred}1.40e-14} &
{\color{SPECgreen!60!black} 1.63} 
\\
968 &rdb968 &{\color{SPECgreen!60!black} 3.034} &{\color{SPECgreen!60!black}9.46e-15} &
{\color{SPECred} 1.95} 
    &{\color{SPECred} 3.657} &{\color{SPECred}1.81e-14} &
{\color{SPECgreen!60!black} 1.72} 
\\
136 &rw136 &{\color{SPECgreen!60!black} 0.021} &{\color{SPECgreen!60!black}4.64e-15} &
{\color{SPECgreen!60!black} 2.76} 
    &{\color{SPECred} 0.024} &{\color{SPECred}8.28e-15} &
{\color{SPECred} 3.15} 
\\
496 &rw496 &{\color{SPECgreen!60!black} 0.718} &{\color{SPECgreen!60!black}8.39e-15} &
{\color{SPECgreen!60!black} 2.72} 
    &{\color{SPECred} 1.248} &{\color{SPECred}1.72e-14} &
{\color{SPECred} 3.06} 
\\
340 &tols340 &{\color{SPECred} 0.109} &{\color{SPECgreen!60!black}4.81e-15} &
{\color{SPECred} 8.43} 
    &{\color{SPECgreen!60!black} 0.078} &{\color{SPECred}6.61e-15} &
{\color{SPECgreen!60!black} 1.95} 
\\
90 &tols90 &{\color{SPECgreen!60!black} 0.003} &{\color{SPECgreen!60!black}3.10e-15} &
{\color{SPECgreen!60!black} 2.01} 
    &{\color{SPECred} 0.003} &{\color{SPECred}4.10e-15} &
{\color{SPECred} 2.24} 
\\
100 &tub100 &{\color{SPECgreen!60!black} 0.005} &{\color{SPECgreen!60!black}3.50e-15} &
{\color{SPECgreen!60!black} 2.34} 
    &{\color{SPECred} 0.006} &{\color{SPECred}5.15e-15} &
{\color{SPECred} 2.82} 
\\
    \bottomrule
  \end{tabular}
\end{center}
\caption{Results special matrices from matrix
  market. \textcolor{SPECgreen!60!black}{Green} entries are better than
  \textcolor{SPECred}{red} entries. }
  \label{tab2}
\end{table}

\end{document}